\documentclass[a4paper, 11pt,russian]{article}
\usepackage[utf8]{inputenc}
\usepackage[english]{babel}
\usepackage{amssymb}
\usepackage{ifthen}
\usepackage{amsmath}

\makeatletter
\def\@seccntformat#1{\csname the#1\endcsname.\ } 
\def\@biblabel#1{#1.} 
\makeatother

\date{}

\def\proofend{$\blacktriangle$\vspace{0.3em}\par}

{\par\addvspace{1mm}{\it Proof\hspace{1.0ex}{#1}.} }%
{\quad$\blacktriangle$\par\addvspace{1mm}}

    \newif\ifNoRemark
    \def\addtheorem#1#2#3#4{ 
    \ifthenelse{\expandafter\isundefined\csname the#2\endcsname}{\newcounter{#2}}{}
    \newenvironment{#1}[1][\global\NoRemarktrue]
     {\par\addvspace{2mm}\noindent 
       \refstepcounter{#2}{\bf #3~\csname the#2\endcsname
      \vphantom{##1}\ifNoRemark.\ \else\ (##1).\fi}\begingroup #4}%
     {\endgroup\par\addvspace{1mm}\global\NoRemarkfalse}
    \expandafter\newcommand\csname b#1\endcsname{\begin{#1}}
    \expandafter\newcommand\csname e#1\endcsname{\end{#1}}
    }

\addtheorem{proposition}{proposition}{Proposition}{\sl}
\addtheorem{theorem}{theorem}{Theorem}{\sl}
\addtheorem{lemma}{lemma}{Lemma}{\sl}
\addtheorem{corol}{corol}{Corollary}{\sl}

\title{On the number of SQS}

\author{Vladimir N. Potapov}

\begin{document}

\maketitle

Abstract. A Steiner quadruple system (briefly $SQS(n)$) is a pair
$(X,B)$ where $|X|=n$ and $B$ is a collection of 4-element blocks
such that every 3-subset of $X$ is contained in exactly one member
of $B$. Hanani \cite{Hanani} proved that  the necessary condition
$n\ {\rm mod}\ 6= 2\ {\rm or}\ 4$ for the existence of a Steiner
quadruple systems of order $n$ is also sufficient. Lenz \cite{Lenz}
proved that the logarithm of the number of different $SQS(n)$ is
greater than $cn^3$ where $c>0$ is a constant and $n$ is admissible.
We prove that the logarithm of the number of different $SQS(n)$ is
$\Theta(n^3\ln n)$ as $n\rightarrow\infty$ and $n\ {\rm mod}\ 6= 2\
{\rm or}\ 4$.


\section{LS and MDS codes}

 By  $Q=[0,q-1]$ denote the subset of integers.  A subset $M$
of $Q^d$ is called an $MDS(t+1,d,q)$ {\it code} (of order
 $q$,
code distance $t+1$ and  length $d$) if $|C\cap\Gamma|=1$ for each
$t$-dimensional face $\Gamma$. These codes achieve equality in the
Singleton bound. As $t=1$, MDS code are equivalent to Latin
$(d-1)$-dimensional cube. If $t=d-2$ then such MDS code is
equivalent to a set of $t$ Mutually  Orthogonal Latin  Squares
(MOLS) of order $q$, and in other cases to a set of $t$ Mutually
Strong Orthogonal Latin $(d-t)$-Cubes. Moreover, a Latin hypercube
is a Cayley table of a multiary qusigroup. A pair of orthogonal
Latin squares corresponds to a pair of orthogonal quasigroups (see
\cite{Pot} or \cite{EMullen}).

By definition MDS code it follows

\begin{proposition}\label{progs1}
Any projection of an MDS code is an MDS code.
\end{proposition}

\begin{proposition}\label{progs2}
Let $M\subset Q^5$ be an MDS code with the code distance $4$ and
$M'$ is a $4$-dimensional projection of $M$. Then there exists an
MDS code $C\subset Q^4$ with  code distance $2$ such that $M'\subset
C$.
\end{proposition}
Proof. By results of \cite{EMullen} any MDS code correspond to a
system of orthogonal quasigroups. So $(x,y,u,v,w)\in M$ whenever
$\qquad \left\{
\begin{array}{l} u=f(x,y);\\
v=g(x,y);\\
w=h(x,y),\\
\end{array} \right.$

\noindent where $f,g,h$ determine a set of $3$ MOLS.

Determine $M'$ by equations $\qquad \left\{
\begin{array}{l} u=f(x,y);\\
v=g(x,y).\\
\end{array} \right.$

Define the function $\varphi:Q^2\rightarrow Q$ by equation
$\varphi(f(x,y),g(x,y))=h(x,y)$.
 The
orthogonality of $f$ and $g$ yields that the function $\varphi$ is
well defined; and the orthogonality  of $f$ and $h$, the
orthogonality of $g$ and $h$ provide that $\varphi$ is a quasigroup.
Hence the set $C=\{(x,y,u,v) \ | \varphi(u,v)=h(x,y)\}$ is an MDS
code and $M'\subset C$ by construction. \proofend

\begin{proposition}\label{progs11}\cite{Wilson}
For every integer $d$ there is an integer  $k(d)$  such that for all
$k>k(d)$ there exists a set of $d$ MOLS of order $k$.

\end{proposition}

Note that  $k(6)$ is not greater  than $75$ \cite{ColDin}.

A subset $T$ of an MDS code $C\subset Q^d$ is called a {\it subcode}
if  $T$ is an MDS code in $A_1\times\dots\times A_d$ and $T=C\cap
(A_1\times\dots\times A_d)$, where $A_i\subset Q$, $i\in
\{1,\dots,d\}$. A definition of a Latin subsquare is analogous.

\begin{proposition}\label{proNMa}
Assume $C$ is an MDS code with a subcode $C_1$ of order $m$, and
assume that a code $C_2$ has the same parameters as $C_1$. Then it
is possible to exchange $C_1$ by $C_2$ in $C$ and to obtain the code
$C'$ with the same parameters as $C$.
\end{proposition}

 A Latin square $f$ is called symmetric  if
$f(x,y)=f(y,x)$ for each $x,y$. It is called nilpotent if $f(x,x)=0$
for every $x$. By using the construction from \cite{Cameron} it is
easy to prove

\begin{proposition}\label{progs12}
Let $q$ be even and $k\leq q/4$. Then there is a symmetric nilpotent
Latin square of order $q$ with subsquare in $K_0\times K_1\times
K_1$ and $K_1\times K_0\times K_1$, where $K_0=[0,q-1]$ and
$K_1=[q-k,q-1]$.
\end{proposition}

\section{Designs}

A {\it $t$-wise balanced design} $t$-BD is a pair $(X,B)$ where $X$
is a finite set of points and $B$ is a set of subsets of $X$, called
blocks, with property that every $t$-element subset of $X$ is
contained in a unique block.
 A {\it $3$-wise bipartite
balanced design} $3$-BBD($n$) is a triple $(X,{g_1,g_2},B)$ where
${g_1,g_2}$ ($|g_1|=|g_2|$) is a partition of $X$, $|X|=n$, $B$ is a
set of $4$-element blocks such that $|b\cap g_i|=2$ for every $b\in
B, i=1,2$ with property that every $3$-element subset $s$ ($s\cap
g\neq \varnothing $) is contained in a unique block.

A {\it Steiner system} $S(t,k,v)$ is a $t$-BD such that $|X|=v$ and
$|b|=k$ for every $b\in B$. If $t=3$ and $k=4$ then this design is
called a Steiner quadruple system. We consider also a $3$-BD denoted
by $S(3,\{4,6\},v)$ consisting of blocks of size $4$ or $6$.

Let $X$ be a set of points, and let $G = \{G_1,\dots,G_d\}$ be a
partition of $X$ into $d$ sets of cardinality $q$. A {\it
transverse} of $G$ is a subset of $X$ meeting each set $G_i$ in at
most
 one point. A set of $w$-element transverses of $G$ is an $H(d,
q, w, t)$ {\it design} (briefly, H-design) if each $t$-element
transverse of $G$ lies in exactly one transverse of the H-design.

An MDS code $M\subset Q^d$ with code distance $t+1$ is equivalent to
$ H(d, q, d, d-t)$, where $G=\{Q_1,\dots, Q_d\}$, $Q_i$ are the
copies of $Q$, and the block $\{x_1,\dots,x_d\}$ lies in the
H-design whenever  $(x_1,\dots,x_d)\in M$.  If $t=2$, an H-design is
called a {\it transversal design}. Transversal designs are
equivalent to systems of MOLS.

If $q$ is even then a $3$-BBD $(X,{g_1,g_2},B)$ is equivalent to the
MDS code $M\subset Q^4$ (with the code distance $2$)  that satisfies
the conditions
\begin{equation}\label{egs1}
 (x,y,u,v)\in M \Rightarrow (y,x,u,v), (x,y,v,u),
(y,x,v,u)\in M;\ \ \forall x,u \in Q\  (x,x,u,u)\in M.
\end{equation}
Here $g_1=Q_1\cup Q_2$, $g_2=Q_3\cup Q_4$, $Q_i$ are copies of $Q$,
and $\{x_1,x_2,x_3,x_4\}\in B$ if $(x_1,\dots,x_4)\in M$ and
$x_1\neq x_2$.


\begin{proposition}\label{thgs1}\cite{Pot}
The logarithm of the  number of MDS codes $M\subset Q^d$  with code
distance $2$ is\footnote{\, Notation $f(x)=\Theta(g(x))$ as
$x\rightarrow x_0$ means that there exist constants $c_2\geq c_1>0$
and a  neighborhood $U$ of $x_0$ such that  for all $x\in U$
$c_1g(x)\leq f(x)\leq c_2g(x)$.} $\Theta(|Q|^{d-1}\ln |Q|)$ as
$n\rightarrow\infty$.
\end{proposition}

Using methods of \cite{Cameron}, \cite{KPS} and  Proposition
\ref{thgs1} we can prove the following theorem.

\begin{theorem}\label{thgs2}
The logarithm of the  number of $3$-wise bipartite balanced designs
on $n$-element set is $\Theta(n^3\ln n)$ as $n\rightarrow\infty$.
\end{theorem}
Proof. Suppose the quasigroup $f$ satisfies the hypothesis  of
Proposition \ref{progs12}. Consider the MDS code $M=\{ (x,y,u,v) \
|\ f(x,y)=f(u,v)\}$. It is easy to see that $M$ meets the conditions
(\ref{egs1}). Furthermore, $M$ has subcodes $B_\sigma$ on
$K_{\sigma_1}\times K_{\sigma_2}\times K_{\sigma_3}\times
K_{\sigma_4}$, where $\sigma=0101, 1001, 0110$ or $1010$.

 For any MDS code $C$
 and permutation $\tau$ we define $C_\tau=
\{(x_{\tau1},\dots,x_{\tau n}) \ |\ x\in C\}$. Let $\Upsilon$ be a
group of permutaions on 4 elements generated by transpositions
$(01)$ and $(23)$.

By Proposition \ref{proNMa} the set $M'=(M\setminus \bigcup_{\tau
\in \Upsilon} K_{\tau(0101)}) \bigcup_{\tau \in \Upsilon} C_\tau $
is an MDS code. By construction, $M'$ satisfies (\ref{egs1}). Since
we use an arbitrary code $C$ of order $k$, the number of $3$-wise
bipartite balanced design is greater than the number of MDS codes of
order $k$. \proofend

The following doubling construction of block designs  is well known
(see \cite{Hartman}).

\begin{proposition}\label{progs3}

\noindent 1. If $S_n\in S(3,4,n)$, $B_n\in 3{\rm -BBD}(n)$ then
there exists $S_{2n}\in S(3,4,2n)$ such that $S_n,B_n\subset
S_{2n}$.

\noindent 2. If $S_n\in S(3,\{4,6\},n)$, $B_n\in 3{\rm -BBD}(n)$
then there exists $S_{2n}\in S(3,\{4,6\},2n)$ such that
$S_n,B_n\subset S_{2n}$.

\end{proposition}

\begin{proposition}\label{progs4}(\cite{Hanani63},  \cite{Hartman} Th.
4.1) There is an injection from $S(3,\{4,6\},n)$ to
$S(3,\{4,6\},2n-2)$.
\end{proposition}

\section{Main results}

The following theorem provides a new construction of SQS based on
MDS codes. Existence of  suitable MDS codes follows from
Propositions \ref{progs1} -- \ref{progs11}.

\begin{theorem}\label{thgs3}

1. If $S_{2n+2}\in S(3,4,2n+2)$, $B_n\in 3{\rm -BBD}(n)$, $n>75$ is
even, then there exists $S_{8n+2}\in S(3,4,8n+2)$ such that
$S_{2n+2},B_n\subset S_{8n+2}$.

2. If $S_{2n+2}\in S(3,\{4,6\},2n+2)$, $B_n\in 3{\rm -BBD}(n)$,
$n>75$ is even, then there exists $S_{8n+2}\in S(3,\{4,6\},8n+2)$
such that $S_{2n+2},B_n\subset S_{8n+2}$.

\end{theorem}

Proof. Below we describe a construction of $S_{8n+2}$ for item 1.
Item 2 is similar.

 Let $I=\{(i,\delta) \ |\ i\in\{0,1,2,3\}, \delta\in
\{0,1\}\}$.
 Denote by $S_8$ a SQS on $I$. Let $S_{10}$ be a SQS on
 $I\cup\{e_1,e_2\}$ such that $\{(i,0),(i,1),e_1,e_2\}\in S_{10}$
 for every $i\in\{0,1,2,3\}$.
Since $n>75$, there exists an $MDS(7,8,n)$ code $M$. We enumerate
these $8$ coordinates by elements of $I$. Consider
$s=\{s_1,s_2,s_3,s_4\}\in S_8$. Denote by $M_s$ the projection of
$M$ on the coordinates $s$. By Proposition~\ref{progs1} $M_s\in
MDS(3,4,n)$. By Proposition~\ref{progs2}, there exists $C_s\in
MDS(2,4,n)$ such that $M_s\subset C_s$.

Now we will construct SQS on a set $\Omega$ where $|\Omega|=8n+2$,
$\Omega=\{e_1,e_2\}\bigcup\limits_{(i,\delta)\in I} A_{(i,\delta)}$
and $|A_{(i,\delta)}|=n$.

Consider H-designs $M^*$, $M_s^*$ and $C_s^*$ with groups
$A_{(i,\delta)}$ that correspond to MDS codes $M$, $M_s$ and $C_s$.
Let us determine quadruples  of four types.

(1) Denote $R_1=\bigcup\limits_{s\in S_8}(C^*_s\setminus M^*_s)$. It
is clear that the blocks of $\bigcup\limits_{s\in S_8}C^*_s$ cover
only once all 3-subsets of $\Omega\setminus \{e_1,e_2\}$ where three
elements lie in different groups.  Besides, a 3-subset  is covered
by a  block of $\bigcup\limits_{s\in S_8} M^*_s$ iff it is included
in a  8-element subset from $M^*$. Note that $\bigcup\limits_{s\in
S_8}(C^*_s)$ and $\bigcup\limits_{s\in S_8}(M^*_s)$ is  H-designs of
type $H(8,n,4,3)$ and $H(8,n,4,2)$, respectively, on
$\Omega\setminus \{e_1,e_2\}$.

(2) Consider any 8-subset $b=\{a^{i,\delta}\in A_{(i,\delta)} \ |
{i,\delta}\in I\}\in M^*$. For every $b\in M^*$ determine a set
$P_b$
 consisting of  blocks $\{a^{s_1},a^{s_2}, a^{s_3}, a^{s_4}\}$, where
$\{s_1, s_2, s_3, s_4\}\in S_{10}$ and blocks $\{a^{s_1},a^{s_2},
a^{s_3}, e_\delta\}$, where $\{s_1, s_2, s_3, \delta\}\in S_{10}$.
Denote by $R_2=\{P_b\ |\ b\in M^*\}$ the set of all these blocks. By
definition of $S_{10}$, the blocks of $R_2$ cover all 3-sets
consisting of $e_1$ or $e_2$ (but not both) and two elements from
$A_{(i,\delta)}$ and $A_{(i',\delta')}$ where $i\neq i'$. Moreover
the blocks of $R_1\cup R_2$ cover all 3-subsets of $\Omega\setminus
\{e_1,e_2\}$, where the three elements lie in different groups.


(3) For any pair $s_0=(i_0,\delta_0)$,  $s_1=(i_1,\delta_1)$ where
$i_0\neq i_1$ consider a $3$-BBD $B_{s_0,s_1}$ with groups $A_{s_0}$
and $A_{s_1}$. Denote $R_3=\bigcup B_{s_0,s_1}$. It is clear that a
3-subset is cover by a block of $R_3$   iff  two elements of the
3-subset lie in $A_{(i,\delta)}$ and the third element lies in
$A_{(i',\delta')}$, where $i\neq i'$.

(4) For $i=0,1,2,3$ consider a Steiner quadruple systems $D_i$ on
the sets $A_{(i,0)}\cup A_{(i,1)}\cup \{e_1,e_2\}$. Define
$R_4=\bigcup D_i$.

By the construction, the blocks from $S_{8n+2}=R_1\cup R_2\cup R_3
\cup R_4$ cover any 3-subset of $\Omega$ only once. To prove
$S_{8n+2}\in S(3,4,8n+2)$, we calculate $|S_{8n+2}|$. It is well
known that SQS of order $m$ consists of $\frac{m(m-1)(m-2)}{4!}$
blocks. Therefore $|R_1|=|S_8|(n^3-n^2)=14(n^3-n^2)$, $R_2=
(|S_{10}|-4)n^2=26n^2$, $R_3=({8 \choose 2}-4)({n \choose
2}n/2)=6n^2(n-1)$, $R_4=4|S_{2n+2}|=(2n+2)(2n+1)n/3$. Then
$$|S_{8n+2}|=|R_1|+|R_2|+|R_3|+|R_4|= 20n^3+6n^2+(2n+2)(2n+1)n/3=$$
$$= 64n^3/3 + 8n^2 + 2n/3= (8n+2)(8n+1)8n/24.$$
\proofend

Note that it is possible to use  SQSs of order $6k+2$ and $6k+4$,
$k\geq 1$ instead of $S_8$ and $S_{10}$.

Now we obtain a lower estimate of the number of block designs as a
corollary of  Propositions \ref{progs3}(2), \ref{progs4}, Theorem
\ref{thgs3}(2) and  the asymptotic estimate from Theorem
\ref{thgs2}.

\begin{theorem}\label{thgs4}
The logarithm of the  cardinality of $S(3,\{4,6\},2n)$  is greater
than $c(n^3\ln n)$, where $c>0$ is a constant.
\end{theorem}
Proof.    If $n$ is  even  then  the statement follows from
Propositions \ref{progs3}(2) and Theorem \ref{thgs2}.

If $n$  is  odd  then we will consider some cases.
 Let $2n=16k+6$. Since $16k+6=2(8k+4)-2$  the
statement follows from Proposition \ref{progs4} and the case of even
$n$. The cases $2n=16k+10=2(2(4k+4)-2)-2$ and $2n=16k+14=2(8k+8)-2$
are simular. If $2n=16k+2$ then we  use Theorems  \ref{thgs2}  and
\ref{thgs3}(2). \proofend

We need some constructions of SQS.

\begin{proposition}\label{progs5}(\cite{Hartman} Th. 4.2)
 There is an injection from $S(3,\{4,6\},n)$ to  $S(3,4,3n-2)$.

\end{proposition}

\begin{proposition}\label{progs6}

1. There is an injection from $S(3,4,n)$ to  $S(3,4,6n-10)$.
(\cite{Hartman} Th. 4.11)

2. If $n\equiv 10\mod{ 12}$ then there exists an injection from
$S(3,4,n)$ to $S(3,4,3n-4)$. (\cite{Hanani} 3.4)
\end{proposition}

The  asymptotic estimate of the number of SQSs is a corollary of
constructions of SQS provided by Propositions \ref{progs3}(1),
\ref{progs5}, \ref{progs6}, Theorem \ref{thgs3}(1) and the
asymptotic estimates from Theorems \ref{thgs2}, \ref{thgs4}.

\begin{theorem}\label{thgs5}
The logarithm of the  cardinality of $S(3,4,n)$  is $\Theta(n^3\ln
n)$ as $n\rightarrow\infty$ and $n\equiv 2\mod{6}$ or $n\equiv
4\mod{6}$.
\end{theorem}

Proof. The upper bound is oblivious (see \cite{Lenz}). To prove
lower bound we will consider apart some subsequences of integers.

(a) Consider a subsequence $n=4k$. For this subsequence the required
asymptotic estimate is a corollary of Theorem \ref{thgs2} and
Proposition \ref{progs3}(1).

(b) Consider the subsequence $n\equiv4(\mod{ 6})$. Then
$n=3(2t+2)-2$ and the required asymptotic estimate is a corollary of
Theorem \ref{thgs4} and Proposition \ref{progs5}.

It retains to  consider three subsequences $n \mod{ 36} = 2, 14\
{\rm or}\ 26$.

(c) If $n=3(12t+10)-4$ then for establishing the required asymptotic
estimate we  use Proposition \ref{progs6}(1) and the proved case
(b).

(d) If $n=6(6t+4)-10$ then  we  use Proposition \ref{progs6}(2) and
the proved case (b).

(e) Consider the case $n\mod{36} = 2$. If $n=6^4t+2= 8(3^42t)+2$
then the required  asymptotic estimate is a corollary of Theorems
\ref{thgs2} and \ref{thgs3}(1). The other cases are reduced to the
subsequence $n=6^4t+2$ by applying Proposition \ref{progs6}(2).
\proofend

\end{document}